\documentclass{article}%
\usepackage{amsmath}%
\usepackage{amsfonts}%
\usepackage{amssymb}%
\usepackage{graphicx}

\begin{document}

\title{On Solving the Cauchy Problem with Propagators}
\author{Henrik Stenlund\thanks{The author is grateful to Visilab Signal Technologies for supporting this work.}
\\Visilab Signal Technologies Oy, Finland}
\date{November 5, 2014}
\maketitle

\begin{abstract}
The abstract first order Cauchy problem is solved in terms of Taylor's series leading to a series of operators which is a propagator. It is found that higher order Cauchy problems can be solved in the same way. Since derivatives of order lower than the Cauchy problem are built-in to the problem, it suffices to solve the higher order derivatives in terms of the lower order ones, in a simple way. \footnote{Visilab Report \#2014-11}
\subsection{Keywords}
Cauchy problem, differential equations, initial value problems
\subsection{Mathematical Classification}
Mathematics Subject Classification 2010: 47D09, 34A12, 35E15, 35F10, 35F40, 35G10, 35G40
\end{abstract}
\subsection{}
Dedicated to my wife Anna-Maija.
\tableofcontents
\section{Introduction}
The first order Cauchy problem (CP) is the following differential equation, treated as an initial value problem. This is a general abstract case with no particular application in mind with ($t \in{R^+},x \in{R}$). 
\begin{equation}
\frac{\partial{u(x,t)}}{\partial{t}}=L{u(x,t)}  \label{eqn10}
\end{equation}
Here $L$ is a linear operator and independent of $t$ and $\frac{\partial}{\partial{t}}$. It is a function of $x$ and the differential operator $\frac{\partial}{\partial{x}}$ but can adopt many forms depending on the particular application. 
\begin{equation}
L=L(x,\frac{\partial}{\partial{x}})  \label{eqn20}
\end{equation}
The initial value of the problem is 
\begin{equation}
u(x,t_0)=u_0(x). 
\end{equation}
The exponential function of $L$ is defined as a Taylor's series with $a \in {C}$, a parameter
\begin{equation}
e^{aL}=\sum^{\infty}_{n=0}\frac{a^n{L^n}}{n!}  \label{eqn30}
\end{equation}
and $L$ commutes with it. The traditional way of solving the (CP) formally is to integrate equation (\ref{eqn10}), see e.g. \cite{Pechukas1966}, leading to
\begin{equation}
u(x,t)=e^{(t-t_o)L}{u_0}, \ \ \ t > t_0  \label{eqn50}
\end{equation}
It is not difficult to prove that this is a unique solution to the (CP). Many functions satisfy equation (\ref{eqn10}). For example, those do, formed from a solution $u$ 
\begin{equation}
\frac{\partial{u(x,t)}}{\partial{t}}  \label{eqn60}
\end{equation}
and
\begin{equation}
h(L){u(x,t)}  \label{eqn70}
\end{equation}
where $h()$ can be almost any function of $L$. Those functions do not, however, satisfy the initial value condition and therefore are not solutions to the (CP).

The (CP) can be treated in multiple dimensions and complex variables too. The (CP) is studied widely and the literature is extensive. The methodologies generally applied are various semigroup methods, regularization methods and abstract distribution methods. The Laplace transform can be applied successfully in many cases. 

The motivation for this paper is to find a solution method for arbitrary order (CP) for solving physical problems. In the following chapter we present the first order (CP) solved with the new approach. Then we apply the method to the second order, to the third, fourth and arbitrary orders. Chapter 3 shows a simple example. We are working in an intuitive manner and are not stepping into rigorous formalism. We keep the presentation simple.

\section{Propagator Solutions}
\subsection{The First Order Cauchy Problem}
We assume well-posedness for the (CP) but do not intend to handle that issue since it has already been treated thoroughly, for example \cite{Pechukas1966}, \cite{Fattorini1983}, \cite{Melnikova2001} and \cite{Xiao1998} and references therein. We presume that the (CP) solution $u(x,t)$ can be expanded as a Taylor's power series, ($t,x \in{R}$) as
\begin{equation}
u(x,t)=\sum^{\infty}_{n=0}{\frac{(t-t_0)^{n}}{n!}[\frac{\partial{^{n}u(x,t)}}{\partial{t^n}}]_{t=t_0}} \label{eqn1000}
\end{equation}
We could equally work in ($t,x \in{C}$) having the requirement of analyticity for the solution. At the end we could flatten our solution to ($t,x \in{R}$). For real variables the requirements are within the interval in question:

\itshape
- $u(x,t)$ is continuous as a function of $t, t > t_0$ 

- the series converges

- all the derivatives exist for any $n=1,2,3..$:
\begin{equation}
\left|\frac{\partial{^{n}u(x,t)}}{\partial{t^n}}\right| < \infty  \nonumber
\end{equation}
\normalfont

The initial value function is
\begin{equation}
u(x,t_0)=u_0(x).
\end{equation}
In order to solve the (CP) we can integrate equation (\ref{eqn10}) and place the expression (\ref{eqn1000}) to it, obtaining 
\begin{equation}
u(x,t)=u_0+L\int^{t}_{t_0}{ds\cdot{u(x,s)}}  \label{eqn1050}
\end{equation}
\begin{equation}
=u_0+L\sum^{\infty}_{n=0}{\frac{1}{n!}[\frac{\partial{^{n}u(x,t)}}{\partial{t^n}}]_{t=t_0}}\int^{t}_{t_0}{ds(s-t_0)^{n}}  \label{eqn1070}
\end{equation}
\begin{equation}
=u_0+L\sum^{\infty}_{n=0}{\frac{(t-t_0)^{n+1}}{(n+1)!}[\frac{\partial{^{n}u(x,t)}}{\partial{t^n}}]_{t=t_0}}  \label{eqn1090}
\end{equation}
Next we need to solve for all the derivatives of the function $u(x,t)$. In effect, we are able to do it by repeatedly differentiating equation (\ref{eqn10}) getting
\begin{equation}
\frac{\partial^n{u(x,t)}}{\partial{t^n}}=L^n{u(x,t)}  \label{eqn1110}
\end{equation}
and thus we get after substitution
\begin{equation}
u(x,t)=u_0+\sum^{\infty}_{n=0}{\frac{L^{n+1}(t-t_0)^{n+1}}{(n+1)!}u_0}=e^{(t-t_0)L}u_0, \ \ \ t > t_0   \label{eqn1140}
\end{equation}
This is the traditional result. The operator $e^{(t-t_0)L}$ has become a propagator for this particular (CP) with $L$. 

\subsection{The Second Order Cauchy Problem}
The second order (CP) is the following initial value problem
\begin{equation}
\frac{\partial^2{u(x,t)}}{\partial{t^2}}=M{u(x,t)}  \label{eqn1210}
\end{equation}
with
\begin{equation}
u(x,t_0)=u_0(x)  \ \ \ \label{eqn1214}
\end{equation}
and 
\begin{equation}
[\frac{\partial{u(x,t)}}{\partial{t}}]_{t=t_0}=u_1(x) \ \ \ \label{eqn1216} 
\end{equation}
being the initial values. $M$ is a linear operator and not a function of $t$ nor $\frac{\partial}{\partial{t}}$. 
\begin{equation}
M=M(x,\frac{\partial}{\partial{x}})  \label{eqn1220}
\end{equation}
We attempt to generalize the procedure of the preceding section. We generate the derivatives in the same way as above, obtaining
\begin{equation}
\frac{\partial{^{n}u(x,t)}}{\partial{t^n}}=   \nonumber
\end{equation}
\begin{equation}
{M^{\frac{n}{2}}u(x,t), n \ even}   \nonumber
\end{equation}
\begin{equation}
{M^{\frac{n-1}{2}}\frac{\partial{u(x,t)}}{\partial{t}}, n \ odd}  \ \ \ \label{eqn1240}
\end{equation}
By applying the series in equation (\ref{eqn1000}) and integrating twice we get
\begin{equation}
u(x,t)=u_0+(t-t_0)u_1+M\sum^{\infty}_{n=0}{\frac{(t-t_0)^{n+2}}{(n+2)!}[\frac{\partial{^{n}u(x,t)}}{\partial{t^n}}]_{t=t_0}}  \label{eqn1250}
\end{equation}
By placing the derivatives and after changing the indices, we have
\newpage

\begin{equation}
u(x,t)=u_0+(t-t_0)u_1+\sum^{\infty}_{n=1}{{(\sqrt{M})}^{2n}\frac{(t-t_0)^{2n}}{(2n)!}u_0}  \nonumber
\end{equation}
\begin{equation}
+\frac{1}{\sqrt{M}}\sum^{\infty}_{n=1}{{(\sqrt{M})}^{2n+1}\frac{(t-t_0)^{2n+1}}{(2n+1)!}u_1}  \label{eqn1260}
\end{equation}
We condense the expressions to get
\begin{equation}
u(x,t)=\sum^{\infty}_{n=0}{\frac{{(\sqrt{M})}^{2n}(t-t_0)^{2n}}{(2n)!}u_0}+\frac{1}{\sqrt{M}}\sum^{\infty}_{n=0}{\frac{{(\sqrt{M})}^{2n+1}(t-t_0)^{2n+1}}{(2n+1)!}u_1}, \ \ \ t > t_0   \label{eqn1270}
\end{equation}
The series are identified as $sinh()$ and $cosh()$ and thus we have an optional representation for the result
\begin{equation}
u(x,t)=cosh[(t-t_0)\sqrt{M}]u_0+\frac{sinh[(t-t_0)\sqrt{M}]}{\sqrt{M}}u_1, \ \ \ t > t_0   \label{eqn1290}
\end{equation}
We have obtained something which is reminiscent of the integrated sine and cosine functions $S(x)$ and $C(x)$, \cite{Kisynski1972}, \cite{Arendt1989}, \cite{Zheng1996}. 

\subsection{The Third Order Cauchy Problem}
Continuing the generalization process, the third order (CP) is the following initial value problem
\begin{equation}
\frac{\partial^3{u(x,t)}}{\partial{t^3}}=P{u(x,t)} ,  \label{eqn1410}
\end{equation}
\begin{equation}
u(x,t_0)=u_0 
\end{equation}
\begin{equation}
[\frac{\partial{u(x,t)}}{\partial{t}}]_{t=t_0}=u_1
\end{equation}
\begin{equation}
[\frac{\partial^2{u(x,t)}}{\partial{t^2}}]_{t=t_0}=u_2
\end{equation}
as the initial values. 
$P$ is a linear operator independent of $t$ and $\frac{\partial}{\partial{t}}$. 
\begin{equation}
P=P(x,\frac{\partial}{\partial{x}})  \label{eqn1420}
\end{equation}
The derivatives are
\begin{equation}
\frac{\partial{^{n}u(x,t)}}{\partial{t^n}}= \nonumber
\end{equation}
\begin{equation}
{P^{\frac{n}{3}}u(x,t), \ \ \ n=0,3,6,9...}  \nonumber
\end{equation}
\begin{equation}
{P^{\frac{n-1}{3}}\frac{\partial{u(x,t)}}{\partial{t}}, \ \ \ n=1,4,7,10...}  \nonumber
\end{equation}
\begin{equation}
{P^{\frac{n-2}{3}}\frac{\partial{^{2}u(x,t)}}{\partial{t^2}}, \ \ \ n=2,5,8,11...}  \label{eqn1440}
\end{equation}
By applying the series in equation (\ref{eqn1000}) we get
\begin{equation}
u(x,t)=\sum^{\infty}_{n=0}{\frac{((t-t_0)\sqrt[3]{P})^{3n}}{(3n)!}}u_0  \nonumber
\end{equation}
\begin{equation}
+\frac{1}{P^{\frac{1}{3}}}\sum^{\infty}_{n=0}{\frac{((t-t_0)\sqrt[3]{P})^{3n+1}}{(3n+1)!}}u_1  \nonumber
\end{equation}
\begin{equation}
+\frac{1}{P^{\frac{2}{3}}}\sum^{\infty}_{n=0}{\frac{((t-t_0)\sqrt[3]{P})^{3n+2}}{(3n+2)!}}u_2, \ \ \ t > t_0   \label{eqn1460}
\end{equation}
\subsection{The Fourth Order Cauchy Problem}
Going further, the fourth order (CP) is the following initial value problem. 
\begin{equation}
\frac{\partial^4{u(x,t)}}{\partial{t^4}}=K{u(x,t)},  \label{eqn1510}
\end{equation}
\begin{equation}
u(x,t_0)=u_0 
\end{equation}
\begin{equation}
[\frac{\partial{u(x,t)}}{\partial{t}}]_{t=t_0}=u_1
\end{equation}
\begin{equation}
[\frac{\partial^2{u(x,t)}}{\partial{t^2}}]_{t=t_0}=u_2
\end{equation}
\begin{equation}
[\frac{\partial^3{u(x,t)}}{\partial{t^3}}]_{t=t_0}=u_3
\end{equation}
are the initial values. $K$ is a linear operator independent of $t$ and $\frac{\partial}{\partial{t}}$. 
\begin{equation}
K=K(x,\frac{\partial}{\partial{x}})  \label{eqn1520}
\end{equation}
The derivatives are
\begin{equation}
\frac{\partial{^{n}u(x,t)}}{\partial{t^n}}= \nonumber
\end{equation}
\begin{equation}
{K^{\frac{n}{4}}u(x,t), \ \ \ n=0,4,8,12...}  \nonumber
\end{equation}
\begin{equation}
{K^{\frac{n-1}{4}}\frac{\partial{u(x,t)}}{\partial{t}}, \ \ \ n=1,5,9,13...}  \nonumber
\end{equation}
\begin{equation}
{K^{\frac{n-2}{4}}\frac{\partial{^{2}u(x,t)}}{\partial{t^2}}, \ \ \ n=2,6,10,12...}  \nonumber
\end{equation}
\begin{equation}
{K^{\frac{n-3}{4}}\frac{\partial{^{3}u(x,t)}}{\partial{t^3}}, \ \ \ n=3,7,11,15...} \label{eqn1540}
\end{equation}
By using the procedure shown we get
\newpage
\begin{equation}
u(x,t)=\sum^{\infty}_{n=0}{\frac{((t-t_0)\sqrt[4]{K})^{4n}}{(4n)!}}u_0  \nonumber
\end{equation}
\begin{equation}
+\frac{1}{K^{\frac{1}{4}}}\sum^{\infty}_{n=0}{\frac{((t-t_0)\sqrt[4]{K})^{4n+1}}{(4n+1)!}}u_1  \nonumber
\end{equation}
\begin{equation}
+\frac{1}{K^{\frac{2}{4}}}\sum^{\infty}_{n=0}{\frac{((t-t_0)\sqrt[4]{K})^{4n+2}}{(4n+2)!}}u_2  \nonumber
\end{equation}
\begin{equation}
+\frac{1}{K^{\frac{3}{4}}}\sum^{\infty}_{n=0}{\frac{((t-t_0)\sqrt[4]{K})^{4n+3}}{(4n+3)!}}u_3, \ \ \ t > t_0   \label{eqn1560}
\end{equation}
\subsection{The $N$th Order Cauchy Problem}
We are able to solve the abstract $N$th order (CP) below with the same assumptions as above. 
\begin{equation}
\frac{\partial^N{u(x,t)}}{\partial{t^N}}=G{u(x,t)}  \label{eqn2510}
\end{equation}
We mark the initial value functions as follows ($i=0,1,2...N-1$)
\begin{equation}
[\frac{\partial^i{u(x,t)}}{\partial{t^i}}]_{t=t_0}=u_i(x)  \label{eqn2500}
\end{equation}
$G$ is a linear operator independent of $t$ and $\frac{\partial}{\partial{t}}$
\begin{equation}
G=G(x,\frac{\partial}{\partial{x}})  \label{eqn2520}
\end{equation}
The series is
\begin{equation}
u(x,t)=\sum^{\infty}_{n=0}{\frac{(t-t_0)^{n}}{n!}[\frac{\partial{^{n}u(x,t)}}{\partial{t^n}}]_{t=t_0}} \label{eqn2525}
\end{equation}
Working as before we get the derivatives
\begin{equation}
\frac{\partial{^{n}u(x,t)}}{\partial{t^n}}= \nonumber
\end{equation}
\begin{equation}
{G^{\frac{n}{N}}u(x,t), \ \ \ n=0,N,2N,3N...}  \nonumber
\end{equation}
\begin{equation}
{G^{\frac{n-1}{N}}\frac{\partial{u(x,t)}}{\partial{t}}, \ \ \ n=1,N+1,2N+1,3N+1...}  \nonumber
\end{equation}
\begin{equation}
{G^{\frac{n-2}{N}}\frac{\partial{^{2}u(x,t)}}{\partial{t^2}}, \ \ \ n=2,N+2,2N+2,3N+2...}  \nonumber
\end{equation}
...
\begin{equation}
{G^{\frac{n-N+1}{N}}\frac{\partial{^{N-1}u(x,t)}}{\partial{t^{N-1}}}, \ \ \ n=N-1,N+N-1,2N+N-1...} \label{eqn2540}
\end{equation}
The solution will be
\begin{equation}
u(x,t)=\sum^{\infty}_{n=0}{\sum^{N-1}_{j=0}{\frac{((t-t_0)G^{\frac{1}{N}})^{Nn+j}G^{\frac{-j}{N}}}{(Nn+j)!}}}u_j(x)  \label{eqn2560}
\end{equation}
or
\begin{equation}
u(x,t)=\sum^{N-1}_{j=0}{(t-t_0)^j\sum^{\infty}_{n=0}{\frac{((t-t_0)^{N}G)^n}{(Nn+j)!}}}u_j(x), \ \ \ t > t_0   \label{eqn2580}
\end{equation}
Since the problem itself provides the derivatives of order ($0...N-1$) at $t_0$, it suffices to solve any higher order ($n>N-1$) derivatives in terms of lower order derivatives, equation (\ref{eqn2540}). The equation (\ref{eqn2580}) above shows that the fractional power operators disappear in the arbitrary order too and we get integer powers. The author believes this approach is new, working just as well for multidimensional problems. 

\section{An Example}
\subsection{The Classic Wave Equation}
One can show that the following exponential propagator relation holds for $f(x)$ if it can be expanded as a Taylor's power series with $x,a \in {C}$.
\begin{equation}
e^{a\frac{\partial}{\partial{x}}}f(x)=f(x+a)  \label{eqn2010}
\end{equation}
This propagator is a universal translation operator. As an application of equation (\ref{eqn1290}) we solve a simple second-order wave equation-like initial value problem, as in equation (\ref{eqn1210}) with
\begin{equation}
M=v^2{\frac{\partial^2}{\partial{x^2}}}  \label{eqn2030}
\end{equation}
The initial values are
\begin{equation}
u(x,t_0)=f(x)  \label{eqn2032}
\end{equation}
\begin{equation}
[\frac{\partial}{\partial{x}}u(x,t)]_{t=t_0}=g(x)  \label{eqn2034}
\end{equation}
We use the formal solution and place the items and change the hyperbolic functions to exponential functions, getting
\begin{equation}
u(x,t)=\frac{1}{2}(e^{v(t-t_0){\frac{\partial}{\partial{x}}}}+e^{{-v(t-t_0){\frac{\partial}{\partial{x}}}}})f(x)+\frac{1}{2v\frac{\partial}{\partial{x}}}(e^{v(t-t_0){\frac{\partial}{\partial{x}}}}-e^{-v(t-t_0){\frac{\partial}{\partial{x}}}})g(x)  \label{eqn2040}
\end{equation}
We can then apply equation (\ref{eqn2010}) and obtain
\begin{equation}
u(x,t)=\frac{1}{2}[f(x+v(t-t_0))+f(x-v(t-t_0))]  \nonumber
\end{equation}
\begin{equation}
+\frac{1}{2v\frac{\partial}{\partial{x}}}[g(x+v(t-t_0))-g(x-v(t-t_0))]  \label{eqn2060}
\end{equation}
The operator
\begin{equation}
\frac{1}{\frac{\partial}{\partial{x}}}  \label{eqn2070}
\end{equation}
is an integration operator and the result will be
\begin{equation}
u(x,t)=\frac{1}{2}[f(x+v(t-t_0))+f(x-v(t-t_0))]  \nonumber
\end{equation}
\begin{equation}
+\frac{1}{2v}\int_{x_0}^{x}{dx[g(x+v(t-t_0))-g(x-v(t-t_0))]}, \ \ \ t > t_0   \label{eqn2090}
\end{equation}
This is the d'Alembert's formula. It has two waves traveling to opposite directions. The interesting second part of it comes from the assumption of nonzero time derivative of the function. 
\section{Discussion}
We have presented a new approach for solving the abstract first order Cauchy problem (CP) producing the familiar result (\ref{eqn1140}). It is based on expanding the solution as a Taylor's time power series. The series becomes an operator series being a propagator for time development turned to the initial value function. 

The procedure extends itself naturally to higher orders, (\ref{eqn1290}), (\ref{eqn1460}) and (\ref{eqn1560}) and to arbitrary order (\ref{eqn2580}). The series become operator series and are propagators for the initial value functions. Identifying the propagator series as elementary functions is possible in some cases, becoming more difficult at higher orders. If the operator has any fractional derivatives, it may still have a formal solution according to this method. However, not all fractional (CP)'s have a solution, see \cite{Bazhlekova1998} for a treatment.

In an arbitrary order $N$ of (CP) the approach is based on solving the derivatives required in the Taylor's series in terms of lower order ($0...N-1$) derivatives in a simple cyclic pattern, equation (\ref{eqn2540}). The fractional operators in the derivatives disappear from the final expressions. The lower order derivatives are provided by the (CP) itself. The main result of this work is equation (\ref{eqn2580}). 



\end{document}